\documentclass[10pt]{article}
\usepackage{}
\usepackage{amsfonts}
\usepackage{mathrsfs}
\usepackage{amssymb}
\usepackage{amsmath,amsthm,amscd,amssymb,mathrsfs,graphicx,amsfonts,mathrsfs}
\usepackage{amsmath}
\usepackage[all]{xy}

\newtheorem{theorem}{\hspace{0.7cm}Theorem }[section]

\newtheorem{lemma}[theorem]{\hspace{0.7cm}Lemma}
\newtheorem{proposition}[theorem]{\hspace{0.7cm}Proposition }
\newtheorem{definition}[theorem]{\hspace{0.7cm}Definition }

\newtheorem{remark}[theorem]{\hspace{0.7cm}Remark }

\def\Ker{\mathop{\rm Ker}\nolimits}
\def\Coker{\mathop{\rm Coker}\nolimits}

\def\Mod{\mathop{\rm Mod}\nolimits}
\def\Hom{\mathop{\rm Hom}\nolimits}

\def\Con{\mathop{\rm Con}\nolimits}

\def\inf{\mathop{\rm inf}\nolimits}
\def\sup{\mathop{\rm sup}\nolimits}

\def\ppd{\mathop{\rm p.pd}\nolimits}
\def\pid{\mathop{\rm p.id}\nolimits}
\def\Pext{\mathop{\rm Pext}\nolimits}
\def\Id{\mathop{\rm Id}\nolimits}
\def\H{\mathop{\rm H}\nolimits}
\def\colim{\mathop{\rm colim}\nolimits}
\def\Ho{\mathop{\rm Ho}\nolimits}

\def\Pgldim{\mathop{\rm p.gldim}\nolimits}

\title{\large \bf On Pure Derived Categories
}
\author{Yuefei Zheng\thanks{{\it E-mail address}: wucaimengr@163.com}
\ and Zhaoyong Huang\thanks{{\it E-mail address}: huangzy@nju.edu.cn}\\
{\footnotesize \it Department of Mathematics, Nanjing
University, Nanjing 210093, Jiangsu Province, China}}
\date{}
\begin{document}
\baselineskip=16pt \maketitle

\begin{abstract}
We investigate the properties of pure derived
categories of module categories, and show that pure derived categories share many nice
properties of classical derived categories. In particular, we show
that bounded pure derived categories can be realized as certain homotopy categories.
We introduce the pure projective (resp. injective) dimension of complexes in pure derived categories,
and give some criteria for computing
these dimensions in terms of the properties of pure projective (resp.
injective) resolutions and pure derived functors. As a consequence,
we get some equivalent characterizations for the finiteness of the pure global
dimension of rings. Finally, pure projective (resp. injective) resolutions of unbounded complexes are considered.
\end{abstract}

\vspace{0.5cm}

\bigskip
\section{ Introduction}
\setcounter{equation}{0}

\vspace{0.2cm}
Let $(\mathcal{A},\mathcal{E})$ be an exact category in the sense of [Q] and $\mathbf{K}(\mathcal{A})$
its homotopy category. Then one can consider the triangulated quotient of $\mathbf{K}(\mathcal{A})$
by $\mathcal{E}$, called the derived category of $(\mathcal{A},\mathcal{E})$, which was studied by
Neeman in [N1]. Now let $R$ be a ring and $R$-Mod the category of left $R$-modules. It is known that
there are two interesting exact structures in $R$-$\Mod$; one is the usual and the other is the pure
exact structure. The derived category with respect to the first one is traditional which provides
a broader framework for studying homological algebra, and to the second one is the pure derived
category which has attracted many authors, see [CHo], [EGO], [Gi], [Kr], [N3], [St] for the details.

In general, triangulated quotients are not intuitive since they are usually realized as calculus of
fractions. However, bounded derived categories are well understood since they are equivalent to
certain homotopy categories of projective modules. It is known that pure projective modules are
exactly projective objects with respect to the pure exact structure, see [KS], [EJ], [P], and [W].
So, it is expected that bounded pure derived categories will share some nice
properties of classical bounded derived categories. In Section 3, we show that for a ring $R$, $R$-$\Mod$
is a full subcategory of its bounded pure derived category.
Moreover, we show that the bounded pure derived category of $R$-$\Mod$ is triangulated equivalent to a triangulated
full subcategory of the bounded above (resp. below) homotopy category of pure projective (resp. injective)
$R$-modules. Note that the results in this section are standard analogs of the corresponding classical ones.

In Section 4, we devote to building triangulated functors from (bounded) pure derived categories.
A very natural choice is the right ``derived" version of $\Hom$. For this, we first establish the pure
projective (resp. injective) resolutions of bounded complexes, and then use them to define right pure derived
functors of $\Hom$ which preserve the corresponding triangles. As applications,
we introduce and study the pure projective (resp. injective) dimension of complexes. In particular,
we obtain some criteria for computing this dimension
in terms of the properties of pure projective (resp. injective) resolutions and the vanishing
of pure derived functors. As a consequence, we get some equivalent characterizations for the finiteness of
the pure global dimension of rings. The results in this section are standard analogs of main results in [AvF],
and generalize the corresponding ones for modules in [KS] and [S].

In Section 5, pure projective (resp. injective) resolutions of certain unbounded complexes are considered.
We use the technique of homotopy (co)limits to show that any bounded below (resp. above) complex
admits a pure projective (resp. injective) resolution.

\bigskip
\section{Preliminaries }
\setcounter{equation}{0}

Throughout this paper, $R$ is an associate ring with identity and
$R$-$\Mod$ is the category of left $R$-modules.
As usual, we use $\mathbf{C}(R)$ and $\mathbf{K}(R)$ to denote
the category of complexes and homotopy category of $R$-$\Mod$, respectively. When we say ``$R$-module",
without an adjective, we mean left $R$-module. For any $X\in\mathbf{C}(R)$, we write
$$\CD
X:=\ \cdots @> >> X^{i-1} @>d_{X} ^{i-1}>> X^{i} @>d_{X} ^{i}>> X^{i+1}
@>d_{X} ^{i+1}>> X^{i+2} @> >> \cdots.\endCD$$ We regard an
$R$-module $M$ as the stalk complex, that is, a complex concentrated
in degree 0.

We recall the bounded conditions for complexes which are standard in homological algebra,
see for example [GM]. Let $X\in\mathbf{C}(R)$. If $X^{i}=0$ for $i\gg 0$,
then $X$ is called {\it bounded above} (or
{\it bounded on the right}). If $X^{i}=0$ for $i\ll 0$, then $X$
is called {\it bounded below} (or {\it bounded on the left}). $X$ is called
{\it bounded} if it is bounded above and below. A cochain map
$f:X\rightarrow Y$ in $\mathbf{C}(R)$ is called a {\it quasi-isomorphism} if it induces
isomorphic homology groups; and $f$ is called a {\it homotopy
equivalence} if there exists a cochain map $g:Y\rightarrow X$ such
that there exist homotopies  $g\circ f\sim \Id_{X}$ and $f\circ g\sim
\Id_{Y}$. For $\Con (f)$ we mean the {\it mapping cone} of a cochain
map $f$. Let $X,Y\in\mathbf{C}(R)$. We use
$\Hom_{R}(X,Y)$ to denote the {\it total complex}, that is, a complex of
$\mathbb{Z}$-modules (where $\mathbb{Z}$ is the additive group of integers)
$$\CD\cdots @> >>\prod \limits_{i\in \mathbb{Z}}\Hom_{R}(X^{i},Y^{i+n}) @>d ^{n}>>
\prod \limits_{i\in \mathbb{Z}} \Hom_{R}(X^{i},Y^{i+n+1}) @> >>
\cdots,
\endCD$$
where $\prod \limits_{i\in \mathbb{Z}}\Hom_{R}(X^{i},Y^{i+n})$ lies
in degree $n$. For any $\varphi \in \Hom_{R}(X,Y)^{n}$,
$d^{n}(\varphi)=(d _{Y}^{i+n}\circ
\varphi^{i}-(-1)^{n}\varphi^{i+1}\circ d ^{i}_{X})_{i\in
\mathbb{Z}}$. Note that this construction defines a bifunctor
$$\Hom_R(-,-):\mathbf{K}(R)^{op}\times \mathbf{K}(R)\to \mathbf{K}(\mathbb{Z}).$$

\begin{definition} {\rm ([W]) A short exact sequence
$$0 \to A \overset{f} \to B \overset{g} \to C \to 0$$
in $R$-$\Mod$ is called {\it pure exact} if for any right
$R$-module $M$, the induced sequence
$$0 \to M\otimes_{R} A\to M \otimes_{R} B \to M \otimes_{R} C \to 0$$
is exact. In this case, $f$ is called {\it pure monic} and $g$ is
called {\it pure epic}.}
\end{definition}

\begin{remark}{\rm
Using the Cohn's theorem (see [R, Theorem 3.69]), we have that a short exact sequence
$$0 \to A  \to B  \to C \to 0$$ in $R$-$\Mod$ is pure exact if and only if
$$0 \to\Hom_{R}(F,A) \to \Hom_{R}(F,B)  \to \Hom_{R}(F,C) \to 0$$ is exact
for any finitely presented $R$-module $F$.}
\end{remark}

\vspace{0.2cm}

In general, the exactness of a complex of $R$-modules is defined ``pointwise". This definition
provides convenience for understanding bounded derived categories.
Let $(\mathcal{A},\mathcal{E})$ be an exact category in the sense of [Q].
Following [N1], a complex $X$ is called {\it acyclic}
with respect to the exact structure of $\mathcal{A}$ if each differential $d_{X}^{i}$
decomposes as $$X^{i}\twoheadrightarrow D^{i} \rightarrowtail X^{i+1},$$ where the former morphism
is admissible epic and the latter one is admissible monic; furthermore, the sequence
$$D^{i}\rightarrowtail X^{i+1}\twoheadrightarrow D^{i+1}$$ is exact for any $i\in \mathbb{Z}$,
see also [Gi, Section 4.2]. Now it is natural for us to propose the following definition,
which provides convenience for understanding bounded pure derived categories later.

\begin{definition} {\rm Let $X\in\mathbf{C}(R)$ and $n\in \mathbb{Z}$. Then $X$ is called
{\it pure exact at $n$,} if the differentials $d_{X}^{n-1}$ and $d_{X}^{n}$ can decompose as above, and the sequence
$$0 \to K^{n} \to X^{n} \to C^{n-1} \to 0$$ is pure exact,
where $K^{n}=\Ker d_{X}^{n}$ and $C^{n-1}=\Coker d_{X}^{n-1}$. $X$ is called
{\it pure exact} if it is pure exact at $n$ for all $n$.}
\end{definition}

\begin{remark}
{\rm \begin{enumerate}
\item[]
\item[(1)] $X\in\mathbf{C}(R)$ is pure exact if and only if
$M\otimes_{R} X$ is exact for any right $R$-module $M$,
and if and only if $\Hom_{R}(F,X)$ is exact for any finitely
presented $R$-module $F$.

\item[(2)] A direct limit of pure exact complexes is again pure exact, since the tensor functor commutes with direct limits
by [R, Theorem 5.27].

\item[(3)] By definition, pure exact complexes coincide with the exact structure in the sense of Neeman [N1].
\end{enumerate}}
\end{remark}

\begin{definition} {\rm ([W]) A module $M\in R$-$\Mod$ is called {\it pure
projective} (resp. {\it injective}) if it is projective (resp. injective)
with respect to every pure exact complex.}
\end{definition}

Let $\mathcal{PP}$ (resp. $\mathcal{PI}$)
be the class of all pure projective (resp.~injective) $R$-modules.
We use $\mathbf{K}^{-}(\mathcal{PP})$ (resp.
$\mathbf{K}^{+}(\mathcal{PI})$) to denote the bounded above
(resp. below) homotopy category of $\mathcal{PP}$ (resp. $\mathcal{PI}$).

\begin{remark} {\rm
\begin{enumerate}
\item[]
\item[(1)] We write $(-)^+:=\Hom_{\mathbb{Z}}(-,\mathbb{Q}/\mathbb{Z})$, where
$\mathbb{Q}$ is the additive group of rational numbers. By [EJ, Proposition 5.3.7],
we have that $M^+$ is a pure injective left $R$-module for any right
$R$-module $M$. Using the fact that
every $R$-module is a direct limit of finitely presented $R$-modules ([R, Lemma 5.39]),
we have that pure projective modules are nothing but summands of
direct sums of finitely presented modules.

\item[(2)] By (1), it is easy to check that a complex $X$ is pure exact
if and only if $\Hom _{R}(P,X)$ is exact for any $P\in \mathcal
{P}\mathcal {P}$, and if and only if $\Hom _{R}(X,I)$ is exact for any
$I\in \mathcal {P}\mathcal {I}$.
\end{enumerate}}
\end{remark}

We need the following definition.

\begin{definition} {\rm A cochain map $f$: $X \to Y$ in $\mathbf{C}(R)$ is called a {\it pure
quasi-isomorphism} if its mapping cone $\Con(f)$ is a pure exact complex.}
\end{definition}

\begin{remark}
{\rm \begin{enumerate}
\item[]
\item[(1)] A cochain map $f: X\to Y$ in $\mathbf{C}(R)$ is a pure
quasi-isomorphism if and only if $$M\otimes_{R}f:
M\otimes_{R}X \to M \otimes_{R}Y$$ is a
quasi-isomorphism for any right $R$-module $M$.

\item[(2)] By Remark 2.6, a cochain map $f$: $X \to Y$ in $\mathbf{C}(R)$ is a pure quasi-isomorphism if and
only if $$\Hom_{R}(P,f):\Hom_{R}(P,X)\rightarrow \Hom_{R}(P,Y)$$ is a
quasi-isomorphism for any $P\in \mathcal {P}\mathcal {P}$, and if and only if $$\Hom_{R}(f,I):\Hom_{R}(Y,I)\to
\Hom_{R}(X,I)$$ is a quasi-isomorphism for any $I\in \mathcal
{P}\mathcal {I}$.
\end{enumerate}}
\end{remark}

The following result concerning both pure exact complexes and pure quasi-isomorphisms is essentially contained in [CFH].

\begin{lemma}
\begin{enumerate}
\item[]
\item[(1)] Let $X\in\mathbf{C}(R)$. Then $X$ is pure exact if and only if
$\emph{Hom}_{R}(P,X)$ is exact for any $P\in \mathbf{K}^{-}(\mathcal
{P}\mathcal {P})$, and if and only if $\emph{Hom}_{R}(X,I)$ is exact
for any $I\in \mathbf{K}^{+}(\mathcal {P}\mathcal {I})$.

\item[(2)] A cochain map $f$ in $\mathbf{C}(R)$ is a pure quasi-isomorphism if and only if
$\emph{Hom}_{R}(P,f)$ is a quasi-isomorphism for any $P\in
\mathbf{K}^{-}(\mathcal {P}\mathcal {P})$, and if and only if
$\emph{Hom}_{R}(f,I)$ is a quasi-isomorphism for any $I\in
\mathbf{K}^{+}(\mathcal {P}\mathcal {I})$.
\end{enumerate}
\end{lemma}

{\it Proof.} The assertion (1) follows from Remark 2.6(2) and [CFH, Lemmas 2.4 and 2.5],
and the assertion (2) follows from Remark 2.8(2) and [CFH, Propositions 2.6 and 2.7].
\hfill$\square$

\begin{lemma}
\begin{enumerate}
\item[]
\item[(1)] Let $f:X\rightarrow Y$ be a pure quasi-isomorphism in $\mathbf{C}(R)$ with $X,~Y\in
\mathbf{K}^{-}(\mathcal {P}\mathcal {P})$. Then $f$ is a homotopy
equivalence.

\item[(2)] Let $f:X\rightarrow Y$ be a pure quasi-isomorphism in $\mathbf{C}(R)$ with $X,~Y\in
\mathbf{K}^{+}(\mathcal {P}\mathcal {I})$. Then $f$ is a homotopy equivalence.
\end{enumerate}
\end{lemma}

{\it Proof.} (1) Because there exists a quasi-isomorphism
$$\Hom_{R}(Y,f):\Hom_{R}(Y,X)\rightarrow \Hom_{R}(Y,Y)$$ by Lemma 2.9, we have
an isomorphism
$$\H^{0}(\Hom_{R}(Y,f)):\H^{0}(\Hom_{R}(Y,X))\rightarrow
\H^{0}(\Hom_{R}(Y,Y)).$$ One can easily check that there exists a
cochain map $g:Y\rightarrow X$ such that $f\circ g\sim$~Id$_{Y}$.
Similarly, there exists a cochain map $h$ such that $g\circ
h\sim$~Id$_{X}$. As a consequence, we have that $g$ and $f$ are homotopy equivalences.

(2) It is the dual of (1). \hfill$\square$

\begin{lemma}
\begin{enumerate}
\item[]
\item[(1)] Let $Y\rightarrow X$ be a pure quasi-isomorphism in $\mathbf{C}(R)$ with $X\in
\mathbf{K}^{b}(R)$ and $Y\in \mathbf{K}^{+}(R)$. Then there exists a
pure quasi-isomorphism $X'\rightarrow Y$ with $X'\in
\mathbf{K}^{b}(R)$.

\item[(2)] Let $X\rightarrow Y$ be a pure quasi-isomorphism in $\mathbf{C}(R)$ with $X\in
\mathbf{K}^{+}(R)$ and $Y\in \mathbf{K}(R)$. Then there exists a pure
quasi-isomorphism $Y\rightarrow X'$ with $X'\in \mathbf{K}^{+}(R)$.
\end{enumerate}
\end{lemma}

{\it Proof.} (1) We can assume that $Y^{n}=0$ for any $n<0$ and that
H$^{i}\Hom_{R}(P,Y)=0$ for any $P\in \mathcal
{P}\mathcal {P}$ and $i\geq m+1$. We have the following commutative diagram
$$\CD
\cdots @> >> 0  @> >> Y^{0} @> >> \cdots @> >> Y^{m-1} @> d_{Y}^{m-1} >> \Ker d_{Y}^{m} @> >> 0 @> >> \cdots \\
@. @. @V \Id_{Y^{0}} VV @. @V \Id_{Y^{m-1}} VV @V   VV @V   VV @.  \\
\cdots @> >> 0 @> >> Y^{0} @> >> \cdots @> >> Y^{m-1} @> d_{Y}^{m-1} >>
Y^{m} @> >> Y^{m+1} @> >> \cdots.
\endCD$$
Let the upper row be the complex $X'$. Since $\Hom_{R}(P,-)$ preserves kernels, the cochain map is clearly a
pure quasi-isomorphism by Remark 2.8(2).

(2) We can assume that $\H^{i}(M\otimes_{R}Y)=0$ for any right $R$-module $M$
and $i\leq -1$. We have the following commutative diagram
$$\CD
\cdots @> >> Y^{-2} @> >> Y^{-1} @> d_{Y}^{-1}>> Y^{0} @>  >> Y^{1} @> >> \cdots \\
@. @V   VV @V   VV @V \widetilde{d}_{Y}^{-1} VV @V \Id_{Y^{1}} VV @. \\
\cdots @> >> 0 @> >> 0 @> >> \Coker d_{Y}^{-1} @> >> Y^{1} @> >> \cdots.
\endCD$$
Let the lower row be the complex $X'$. Since $M\otimes_{R}-$ preserves
cokernels, the cochain map is clearly a
pure quasi-isomorphism by Remark 2.8(1). \hfill$\square$

\bigskip
\section{ Pure derived categories}
\setcounter{equation}{0}

Put $\mathbf{K}_{\mathcal {P}\mathcal {E}}(R):=\{X\in
\mathbf{K}(R)\mid X$~is~pure~exact$\}$. Notice that pure exact
complexes are closed under homotopy equivalences, so
$\mathbf{K}_{\mathcal {P}\mathcal {E}}(R)$ is well defined. If
$f:X\to Y$ is a cochain map between pure exact complexes, then
$\Con(f)$ is again pure exact. Thus $\mathbf{K}_{\mathcal
{P}\mathcal {E}}(R)$ is a triangulated subcategory of
$\mathbf{K}(R)$. Because pure exact complexes are closed under summands by definition,
$\mathbf{K}_{\mathcal{P}\mathcal {E}}(R)$ is a thick subcategory of $\mathbf{K}(R)$.
Then by the Verdier's correspondence, we get the pure derived category
$$\mathbf{D_{pur}}(R):=\mathbf{K}(R)/\mathbf{K}_{\mathcal {P}\mathcal
{E}}(R).$$ Similarly, we define
$$\mathbf{D^{*}_{pur}}(R):=\mathbf{K^{*}}(R)/\mathbf{K^{*}}_{\mathcal
{P}\mathcal {E}}(R)$$ for $*\in \{+,~-,~b\}$. Note that the pure
derived category coincides with the one given in [N1] and pure exact
complexes here are exactly the exact structure there.

Note that, as usual, a morphism from $X$ to $Y$ in $\mathbf{D_{pur}}(R)$ can be
viewed as a graph (left roof)
$$\xymatrix{
 X & \bullet
\ar@{=>}[l]_{s} \ar[r]^{a} & Y}$$ with $s$ a pure quasi-isomorphism ([GM, Chapters III.2.8 and III.2.9]).
Two roofs $$\xymatrix{
 X & \bullet
\ar@{=>}[l]_{s} \ar[r]^{a} & Y}~~\mathrm{and}~~\xymatrix{
 X & \bullet
\ar@{=>}[l]_{s'} \ar[r]^{a'} & Y}$$ are equivalent if there exists the following commutative diagram
$$\xymatrix{
& \bullet \ar@{=>}[dl]_{s} \ar[dr]^{a}\\
X & \bullet  \ar@{=>}[l]_{g} \ar[r]\ar[d] \ar[u] & Y\\
& \bullet \ar@{=>}[ul]^{s'} \ar[ur]_{a'} }$$ with $g$ a pure quasi-isomorphism.
So, two complexes $X, Y$ are isomorphic in $\mathbf{D_{pur}}(R)$ if there exists a graph
$$\xymatrix{ X & \bullet
\ar@{=>}[l]_{s} \ar@{=>}[r]^{a} & Y}$$ with $s$ and $a$ pure quasi-isomorphisms.
If either $Y\in \mathbf{K}^{+}(\mathcal {P}\mathcal{I})$ or $X\in \mathbf{K}^{-}(\mathcal {P}\mathcal {P})$,
then morphisms in $\Hom _{\mathbf{D_{pur}(R)}}(X,Y)$ are easy enough as showed below.

\begin{proposition}
\begin{enumerate}
\item[]
\item[(1)] Let $X\in \mathbf{K}^{-}(\mathcal {P}\mathcal {P})$ and $Y\in \mathbf{K}(R)$. Then the localization functor
$$\mathbb{F}:\Hom_{\mathbf{K}(R)}(X,Y)\rightarrow
\Hom_{\mathbf{D_{pur}}(R)}(X,Y),~f\mapsto f/ \Id_{X}\ (left\ roof),$$
induces an isomorphism of abelian groups.

\item[(2)] Let $Y\in \mathbf{K}^{+}(\mathcal {P}\mathcal {I})$ and $X\in \mathbf{K}(R)$. Then the localization functor
$$\mathbb{F}:\Hom_{\mathbf{K}(R)}(X,Y)\rightarrow
\Hom_{\mathbf{D_{pur}}(R)}(X,Y),~f\mapsto \Id_{Y}\backslash f\
(right\ roof),$$ induces an isomorphism of abelian groups.
\end{enumerate}
\end{proposition}

{\it Proof.} We only need to prove (1). If $f/\Id_{X}=0=0/\Id_{X}$,
then there exists a pure quasi-isomorphism $g:Z \rightarrow X$ such
that $f\circ g\sim 0$. By the proof of Lemma 2.10, there exists a
pure quasi-isomorphism $h:X\to Z$ such that $g\circ h\sim \Id_{X}$.
So $f\sim 0$. For any $f/s\in \Hom_{\mathbf{D_{pur}(R)}}(X,Y)$,
since $s$ is a pure quasi-isomorphism, again by the proof of Lemma
2.10 there exists a pure quasi-isomorphism $t$ such that $s\circ
t\sim \Id_{X}$. So we have $f/s=(f\circ t)/\Id_{X}$ in
$\mathbf{D_{pur}}(R)$. \hfill$\square$

\vspace{0.2cm}

\begin{proposition} For a ring $R$, we have
\begin{enumerate}
\item[(1)] $\mathbf{D^{b}_{pur}}(R)$ is a full subcategory of
$\mathbf{D^{+}_{pur}}(R)$, and $\mathbf{D^{+}_{pur}}(R)$ is a full
subcategory of $\mathbf{D_{pur}}(R)$.

\item[(2)] $\mathbf{D^{b}_{pur}}(R)$ is a full subcategory of
$\mathbf{D^{-}_{pur}}(R)$, and $\mathbf{D^{-}_{pur}}(R)$ is a full
subcategory of $\mathbf{D_{pur}}(R)$.

\item[(3)] $\mathbf{D^{b}_{pur}}(R)=\mathbf{D^{-}_{pur}}(R)\cap
\mathbf{D^{+}_{pur}}(R)$.
\end{enumerate}
\end{proposition}

{\it Proof.} The assertion (1) is a consequence of [GM, Proposition 3.2.10] and Lemma 2.11,
and the assertion (2) is the dual of (1). The assertion (3) is an immediate consequence of (1) and
(2). \hfill$\square$

\vspace{0.2cm}

\begin{theorem}
For a ring $R$, $R$-$\Mod$ is a full subcategory of
$\mathbf{D^{b}_{pur}}(R)$, that is, the composition of functors
$$R{\text -}\Mod \rightarrow \mathbf{K^{b}}(R)\rightarrow
\mathbf{D^{b}_{pur}}(R)$$ is fully faithful.
\end{theorem}

{\it Proof.} For any $X, Y\in R$-$\Mod$, it suffices to
prove that the morphism $$\mathbb{F}:\Hom_{R}(X,Y)\rightarrow
\Hom_{\mathbf{D_{pur}}(R)}(X,Y)$$ is an isomorphism.

Let $f\in \Hom_{R}(X,Y)$. If $\mathbb{F}(f)=0$, then there exists a
pure quasi-isomorphism $s:Z \rightarrow X$ such that
$f\circ s\sim 0$. So $\H^{0}(f)\circ\H^{0}(s)=0$. Since
$\H^{0}(s)$ is an isomorphism, we have $f=0$.

Let $a/s$ be a morphism in $\Hom_{\mathbf{D_{pur}}(R)}(X,Y)$. Then we have a diagram
$$\xymatrix{
X & Z\ar@{=>}[l]_{s} \ar[r]^{a} & Y},$$ where $s$ is a pure
quasi-isomorphism, and hence a quasi-isomorphism. So $\H^{0}(s)\in
\Hom_{R}(\H^{0}(Z),X)$ is an isomorphism in $R$-$\Mod$ (note that $\H^0(X)=X$). Put
$f:=\H^{0}(a)\circ \H^{0}(s)^{-1} \in \Hom_{R}(X,Y)$. Consider the
truncation $$\CD U:=~~~ \cdots @> >> Z^{-2} @> >> Z^{-1} @>d^{-1}>>
\Ker d_{Z}^{0} @> >> 0
\endCD$$ of $Z$ and the canonical map $i:U\rightarrow Z$. Note that, as in Lemma 2.11, $i$ is a pure quasi-isomorphism.
Then $s\circ i$ is also a pure quasi-isomorphism. From the commutative diagram
$$\CD
U @> i >> Z \\
@V   VV @V s VV  \\
\H^{0}(Z) @> \textrm{H}^{0}(s) >> X,
\endCD$$
we get $f\circ s\circ i=\H^{0}(a)\circ \H^{0}(s)^{-1}\circ s\circ
i=a\circ i$. So the following diagram of complexes
$$\xymatrix{
& Z \ar@{=>}[dl]_{s} \ar[dr]^{a}\\
X & U  \ar@{=>}[l]_{si} \ar[r]^{ai}\ar[d]^{si} \ar[u]^{i} & Y\\
& X \ar@{=>}[ul]^{\textrm{Id}_{X}} \ar[ur]_{f} }$$ is
commutative. It follows that
$\mathbb{F}(f)=f/\Id_{X}=a/s$. \hfill$\square$

\vspace{0.2cm}

For any $X\in\mathbf{C}(R)$, we write

$$~~~~\mathbf{inf_{p}}X:=\inf \{n\in \mathbb{Z} \mid X\ {\rm is\ not\ pure\ exact\ at}\ n\},\ {\rm and}$$

$$\mathbf{sup_{p}}X:=\sup \{n\in \mathbb{Z}\mid X\ {\rm is\ not\ pure\ exact\ at}\ n\}.$$
If $X$ is not pure exact at $n$ for any $n$, then we set $\mathbf{inf_{p}}X=-\infty$ and $\mathbf{sup_{p}}X=\infty$.
If $X$ is pure exact at $n$ for all $n$, that is, $X$ is a pure exact complex, then we set $\mathbf{inf_{p}}X=\infty$
and $\mathbf{sup_{p}}X=-\infty$. We will heavily rely on these two numbers in the remainder
of this paper.

Put $$\mathbf{K^{-,pb}}(\mathcal {P}\mathcal {P}):=\{X\in
\mathbf{K^{-}}(\mathcal {P}\mathcal {P}) \mid \mathbf{inf_{p}}X \ {\rm is\
finite}\},\  {\rm and}$$
$$\mathbf{K^{+,pb}}(\mathcal {P}\mathcal {I}):=\{X\in
\mathbf{K^{+}}(\mathcal {P}\mathcal {I}) \mid \mathbf{sup_{p}}X \ {\rm is\
finite}\}.$$

\begin{proposition} Let $X\in\mathbf{C}(R)$. Then the following hold.
\begin{enumerate}
\item[(1)] $X$ is pure exact in degree $\leq n$
if and only if $M\otimes_{R} X$ is exact in degree $\leq n$ for any right $R$-module $M$.

\item[(2)] $X$ is pure exact in degree $\geq n$
if and only if $\Hom_{R}(P,X)$ is exact in degree $\geq n$ for any $P\in \mathcal {PP}$.

\item[(3)] The numbers $\mathbf{inf_{p}}X$ and $\mathbf{sup_{p}}X$ are well defined for
any $X\in \mathbf{D_{pur}}(R)$, that is, if $X\cong Y$ in $\mathbf{D_{pur}}(R)$,
then $\mathbf{inf_{p}}X=\mathbf{inf_{p}}Y$ and $\mathbf{sup_{p}}X=\mathbf{sup_{p}}Y$.

\item[(4)] $\mathbf{K^{-,pb}}(\mathcal {P}\mathcal {P})$
and $\mathbf{K^{+,pb}}(\mathcal {P}\mathcal {I})$ are triangulated subcategories of
$\mathbf{K^{-}}(\mathcal{P}\mathcal {P})$ and $\mathbf{K^{+}}(\mathcal {P}\mathcal{I})$,
respectively.
\end{enumerate}
\end{proposition}

{\it Proof.} (1) Consider the following commutative diagram (tensor products act on $R$)
$$\xymatrix{\cdots \ar[r] & M\otimes X^{n-1}\ar[rr]^{M\otimes
d_{X}^{n-1}} \ar@{>>}[rd]_<<<<<<<{M\otimes \widetilde{d}_{X}^{n-2}} & &
M\otimes X^{n} \ar@{>>}[rd]_<<<<<<<{M\otimes \widetilde{d}_{X}^{n-1}}
\ar[rr]^{M\otimes d_{X}^{n}} & & M\otimes X^{n+1} \ar[r] & \cdots,
\\ & & M\otimes C^{n-2} \ar[ru]_{M\otimes \iota ^{n}}& & M\otimes C^{n-1}
\ar[ru]_{M\otimes \iota ^{n+1}}& & }$$ where $\widetilde{d}_{X}^{n-2}$
(resp. $\widetilde{d}_{X}^{n-1}$) denotes the cokernel of $d_{X}^{n-2}$
(resp. $d_{X}^{n-1}$) and $ \iota ^{n}$ (resp. $\iota ^{n+1}$)
denotes the kernel of $d_{X}^{n}$ (resp. $d_{X}^{n+1}$). Then the assertion follows standardly.

(2) The proof is similar to that of (1).

(3) We only need to prove the assertion whenever both $\mathbf{inf_{p}}X$
(resp. $\mathbf{sup_{p}}X$) and $\mathbf{inf_{p}}Y$ (resp. $\mathbf{sup_{p}}Y$) are finite.
By Remark 2.8, it is an immediate consequence of (1) and (2).

(4) We only prove that $\mathbf{K^{-,pb}}(\mathcal {P}\mathcal {P})$ is a triangulated subcategory of
$\mathbf{K^{-}}(\mathcal{P}\mathcal {P})$ and the proof of the other assertion is similar.
Observe that $\mathbf{K^{-,pb}}(\mathcal {P}\mathcal {P})$ is closed under shifts. So it suffices
to show that $\mathbf{K^{-,pb}}(\mathcal {P}\mathcal {P})$ is closed under extensions. Let
$$X\to Y\to Z\to X[1]$$ be a triangle
in $\mathbf{K^{-}}(\mathcal{P}\mathcal {P})$ with $X,Z\in \mathbf{K^{-,pb}}(\mathcal {P}\mathcal {P})$. Then
we have a triangle
$$X\otimes_R M\to Y\otimes_R M\to Z\otimes_R M\to X[1]\otimes_R M$$ in $\mathbf{K}(\mathbb{Z})$ for any right $R$-module $M$.
It induces a long exact sequence of homological groups since $\H^0(-)$ is cohomological by [GM, Chapter IV.1.6].
By (1) there exists
$n\in \mathbb{Z}$ such that both $X\otimes_R M$ and $Z\otimes_R M$ are exact in degree $\leq n$,
so $Y\otimes_R M$ is also exact in degree $\leq n$. Thus we have $Y\in\mathbf{K^{-,pb}}(\mathcal {P}\mathcal {P})$.
\hfill$\square$

\begin{proposition}
\begin{enumerate}
\item[]
\item[(1)] There exist a functor $P:\mathbf{K^{b}}(R)\rightarrow
\mathbf{K^{-,pb}}(\mathcal {P}\mathcal {P})$ and a pure
quasi-isomorphism $f_{X}:P_{X}\rightarrow X$ for any $X\in
\mathbf{K^{b}}(R)$, which is functorial in $X$.

\item[(2)] There exist a functor $I:\mathbf{K^{b}}(R)\rightarrow
\mathbf{K^{+,pb}}(\mathcal {P}\mathcal {I})$ and a pure
quasi-isomorphism $g_{X}:X\rightarrow I_{X}$ for any $X\in
\mathbf{K^{b}}(R)$, which is functorial in $X$.
\end{enumerate}
\end{proposition}

{\it Proof.} We first prove that for any $X\in \mathbf{K^{b}}(R)$,
there exists a pure quasi-isomorphism $P_{X}\rightarrow X$ with
$P_{X}\in \mathbf{K^{-,pb}}(\mathcal {P}\mathcal {P})$. We proceed
by induction on the cardinal of the finite set $\mathcal
{W}(X):=\{i\in \mathbb{Z}\mid X^{i}\neq 0  \}$.

If $\mathcal {W}(X)=1$, then the assertion follows from the fact
that every module admits a pure projective precover (see [EJ, Example 8.3.2]).

Now suppose that $\mathcal {W}(X)\geq 2$ with $X^{j}\neq 0$ and
$X^{i}=0$ for any $i<j$. Then we have a distinguished triangle
$$\CD
X_{1} @>u>> X_{2} @> >> X @> >> X_{1}[1]
\endCD$$ in $\mathbf{K^{b}}(R)$, where $X_{1}=X^{j}[-j-1]$ and $X_{2}=X^{>j}$. By the induction hypothesis,
there exist pure quasi-isomorphisms $f_{X_{1}}:P_{X_{1}}\rightarrow
X_{1}$ and $f_{X_{2}}:P_{X_{2}}\rightarrow X_{2}$ with $P_{X_{1}},
P_{X_{2}}\in \mathbf{K^{-,pb}}(\mathcal {P}\mathcal {P})$. Then by
Lemma 2.9, $f_{X_{2}}$ induces an isomorphism
$$\Hom_{\mathbf{K}(R)}
(P_{X_{1}},P_{X_{2}})\cong \Hom_{\mathbf{K}(R)}(P_{X_{1}},X_{2}).$$
So there exists a morphism $f:P_{X_{1}}\rightarrow P_{X_{2}}$, which
is unique up to homotopy, such that $f_{X_{2}}\circ f=u\circ f_{X_{1}}$. We
have the distinguished triangle
$$\CD P_{X_{1}} @> f >> P_{X_{2}} @> >> \Con(f) @> >> P_{X_{1}}[1]
\endCD$$ in $\mathbf{K^{-,pb}}(\mathcal{PP})$. Then there exists a morphism
$f_{X}:\Con(f)\rightarrow X$ such that the following diagram
$$\CD
P_{X_{1}} @> f >> P_{X_{2}} @> >> \Con(f) @> >> P_{X_{1}}[1] \\
@V f_{X_{1}} VV @V f_{X_{2}} VV @V f_{X} VV @V f_{X_{1}}[1] VV  \\
X_{1} @> u >> X_{2} @> >> X @> >> X_{1}[1]
\endCD$$
in $\mathbf{K}(R)$ commutes.
For any $P\in \mathcal {P}\mathcal {P}$, we have the following
commutative diagram $$\CD
\Hom_{R}(P,P_{X_{1}}) @> >> \Hom_{R}(P,P_{X_{2}}) @> >> \Hom_{R}(P,\Con(f)) @> >> \Hom_{R}(P,P_{X_{1}}[1])  \\
@V (f_{X_{1}})_{\ast} VV @V (f_{X_{2}})_{\ast} VV @V (f_{X})_{\ast} VV @V (f_{X_{1}}[1])_{\ast} VV   \\
\Hom_{R}(P,X_{1}) @> >> \Hom_{R}(P,X_{2}) @> >> \Hom_{R}(P,~X) @> >>
\Hom_{R}(P,X_{1}[1])
\endCD$$ in
$\mathbf{K}(\mathbb{Z})$,
where both rows are exact triangles and $(-)_{*}$ denotes the functor $\Hom_{R}(P,-)$. Since both $f_{X_{1}}$
and $f_{X_{2}}$ are pure quasi-isomorphisms, we have that
both $(f_{X_{1}})_{\ast}$ and $(f_{X_{2}})_{\ast}$ are quasi-isomorphisms. Passing to
homology we get that $(f_{X})_{\ast}$ is a quasi-isomorphism, so $f_{X}$ is a
pure quasi-isomorphism by Remark 2.8(2).

Put $P_{X}:=\Con(f)$. Then we have a pure quasi-isomorphism
$f_{X}:P_{X}\rightarrow X$ with $P_{X}\in
\mathbf{K^{-,pb}}(\mathcal {P}\mathcal {P})$. In the following we prove that $f_{X}$
is functorial in $X$.

Let $X,Y\in \mathbf{K}^{b}(R)$. Then we have two pure
quasi-isomorphisms $f_{X}:P_{X}\rightarrow X$ and $f_{Y}:P_{Y}\rightarrow
Y$. These induce an isomorphism
$$\Hom_{\mathbf{K}(R)}(P_{X},P_{Y})\cong
\Hom_{\mathbf{K}(R)}(P_{X},Y).$$ Let $f:X\rightarrow Y$ be a cochain
map. Then there exists a cochain map $f\circ f_{X}:P_{X}\rightarrow
Y$. Using the above isomorphism, we have that there exists a unique cochain map
$f':P_{X}\rightarrow P_{Y}$ such that the following diagram
$$\CD
P_{X} @>f_{X}>> X \\
@V f' VV @V f VV  \\
P_{Y} @>f_{Y}>> Y
\endCD$$
commutes up to homotopy. This completes the proof by putting $Y=X$. \hfill$\square$

(2) It is the dual of (1) just using the fact that every module admits
a pure injective preenvelope by [EJ, Proposition 5.3.9].

\begin{theorem} For a ring $R$, there exist triangle-equivalences as follows.
\begin{enumerate}
\item[(1)] $\mathbf{D^{b}_{pur}}(R)\simeq \mathbf{K^{-,pb}}(\mathcal
{P}\mathcal {P})$.

\item[(2)] $\mathbf{D^{b}_{pur}}(R)\simeq \mathbf{K^{+,pb}}(\mathcal
{P}\mathcal {I})$.
\end{enumerate}
\end{theorem}

{\it Proof.} We only need to prove (1). Let $\mathbb{H}$ be the
composition of the embedding $$\mathbf{K^{-,pb}}(\mathcal {P}\mathcal
{P})\hookrightarrow \mathbf{K^{-}}(R)$$ and the localization functor
$$\mathbb{F}:\mathbf{K^{-}}(R)\rightarrow \mathbf{D^{-}_{pur}}(R).$$
For any $X\in \mathbf{K^{-,pb}}(\mathcal {P}\mathcal {P})$, there
exists $n\in \mathbb{Z}$ such that $\mathbf{inf_{p}}X=n$. So $X$
is pure exact in degree $\leq n-1$ and the following cochain map $f$
is a pure quasi-isomorphism.
$$\CD
  ~~~~~X:=~\cdots @> >> X^{n-1} @> >> X^{n} @> >> X^{n+1} @> >> X^{n+2} @> >> \cdots \\
  @V f VV @V   VV @V   VV @V  VV @V  VV @.  \\
  X^{\supset n}:=~~\cdots @> >> 0 @> >> \Coker d^{n-1} @> >> X^{n+1} @> >> X^{n+2} @> >> \cdots.
\endCD$$
It follows that $\mathbb{H}(X)\cong X^{\supset n}$ in
$\mathbf{D_{pur}}(R)$. So $\mathbb{H}(X)\in
\mathbf{D^{b}_{pur}}(R)$ and hence $\mathbb{H}$ induces a functor
from $\mathbf{K^{-,pb}}(\mathcal {P}\mathcal {P})$ to
$\mathbf{D^{b}_{pur}}(R)$, again denoted by $\mathbb{H}$. By
Propositions 3.1 and 3.5, $\mathbb{H}$ is fully faithful
and dense. This completes the proof. \hfill$\square$

\bigskip
\section{Derived functors and dimensions}
\setcounter{equation}{0}

In this section, we introduce and investigate the pure projective and
injective dimensions of complexes based on pure derived functors of $\Hom$ in $\mathbf{D^{b}_{pur}}(R)$.
For the pure projective and injective dimensions of modules and pure derived functors in
$R$-$\Mod$, we refer to [KS] and [S].

We have already known that $\Hom_{R}(P,-)$ transforms pure
quasi-isomorphisms to quasi-isomorphisms for any $P\in
\mathbf{K^{-}}(\mathcal {P}\mathcal {P})$. In order to define pure
projective (resp. injective) resolutions of complexes in
$\mathbf{D_{pur}}(R)$, we need the following lemma.

\begin{lemma}
Let $X$ be a pure exact complex of $R$-modules. Then we have
\begin{enumerate}
\item[(1)] $M\otimes_{R} X$ is a pure exact complex for any
right $R$-module $M$.

\item[(2)] $\Hom_{R}(P,X)$ is a pure exact complex for any
$P\in\mathcal{PP}$.

\item[(3)] $\Hom_{R}(X,I)$ is a pure exact complex for any $I\in
\mathcal{P}\mathcal{I}$.
\end{enumerate}
\end{lemma}

{\it Proof.} (1) It is obvious by the associativity of tensor products.

(2) We will prove that $\Hom_{\mathbb{Z}}(F,\Hom_{R}(P,X))$ is exact
for any finitely presented $\mathbb{Z}$-module $F$ and $P\in\mathcal{PP}$.

By Remark 2.6, $\mathcal{P}\mathcal{P}$ consists of summands of
direct sums of finitely presented $R$-modules. So we may assume that $P$ is
finitely presented. Note that
$P\otimes_{\mathbb{Z}} F$ is a finitely presented
$R$-module. So by the adjoint isomorphism
$\Hom_{\mathbb{Z}}(F,\Hom_{R}(P,X))\cong
\Hom_{R}(P\otimes_{\mathbb{Z}} F,X)$, we have that
$\Hom_{R}(P,X)$ is pure exact.

(3) Let $I\in\mathcal{PI}$. Then $I$ is a direct summand
of $I^{++}$ by [EJ, Proposition 5.3.9]. We may assume $I=M^{+}$ for some right $R$-module $M$.
Let $F$ be a finitely presented $\mathbb{Z}$-module. By the adjoint isomorphism theorem, we have the
isomorphisms
$$\Hom_{\mathbb{Z}}(F,\Hom_{R}(X,M^{+})) \cong
\Hom_{\mathbb{Z}}(F,(M\otimes_{R}X)^+) \cong
(F\otimes_{\mathbb{Z}}M\otimes_{R}X)^+.$$
By (1), $F\otimes_{\mathbb{Z}} M\otimes_{R} X$ is pure
exact. So $(F\otimes_{\mathbb{Z}}M\otimes_{R}X)^+$ is exact, and hence
$\Hom_{R}(X,I)$ is pure exact.
\hfill$\square$

\begin{remark}{\rm By [CFH, Lemmas 2.4 and 2.5] and Lemma 4.1, after a standard computation we have
\begin{enumerate}
\item[(1)] $\Hom_{R}(P,-)$ preserves pure exact complexes for any $P\in
\mathbf{K^{-}}(\mathcal {P}\mathcal {P})$.

\item[(2)] $\Hom_{R}(-,I)$ preserves pure exact complexes for any $I\in
\mathbf{K^{+}}(\mathcal {P}\mathcal {I})$.
\end{enumerate}}
\end{remark}

\begin{definition}{\rm Let $X\in\mathbf{D_{pur}}(R)$.
\begin{enumerate}
\item[(1)] A {\it pure projective resolution} of $X$ is a pure
quasi-isomorphism $f:P \to X$ with $P$ a complex of pure
projective $R$-modules, such that $\Hom_{R}(P,-)$ preserves pure exact
complexes. Dually, a {\it pure injective resolution} of $X$ is
defined.

\item[(2)] $X$ is said to have {\it pure projective dimension} at most $n$, written $\ppd_{R}X\leq n$, if there exists a pure
projective resolution $P\rightarrow X$ with $P^{i}=0$ for any $i<-n$. If $\ppd _{R}X\leq n$ for all $n$, then we write
$\ppd _{R}X=-\infty$; and if there exists no $n$ such that
$\ppd _{R}X\leq n$, then we write $\ppd _{R}X=\infty$. Dually, the
{\it pure injective dimension} $\pid_{R}X$ of $X$ is defined.
\end{enumerate}}
\end{definition}

\begin{remark}
{\rm \begin{enumerate}
\item[]
\item[(1)] Let $X$ be an $R$-module (viewed as a complex concentrated in degree 0),
then these definitions coincide with the usual ones, see [KS] and [S].

\item[(2)] In the above definition, $\ppd _{R}X=-\infty$ means that $X$ is a pure exact complex.
\end{enumerate}}
\end{remark}

\vspace{0.2cm}

These dimensions can be also expressed by the following equalities.
$$\ppd _{R}X=-\sup\{\inf\{n\in \mathbb{Z}\mid P^{n}\neq 0\}\mid P\rightarrow X\
{\rm is\ a\ pure\ projective\ resolution}\}, {\rm and}$$
$$\pid _{R}X=\inf\{\sup\{n\in \mathbb{Z}\mid I^{n}\neq 0\}\mid X\rightarrow I\
{\rm is\ a\ pure\ injective\ resolution}\}.$$

Let $X\in\mathbf{D_{pur}^{b}}(R)$. Then by Proposition 3.5, there exists
a complex $P\in\mathbf{K^{-}}(\mathcal{PP})$ such that $P\cong X$ in $\mathbf{D_{pur}^{b}}(R)$.
By Remark 4.2, $\Hom_{R}(P,-)$
preserves pure exact complexes, and hence preserves pure
quasi-isomorphisms. Then after an easy computation we get a pure
quasi-isomorphism from $P$ to $X$. The statements for the pure injective
version are dual. Thus, if $X\in \mathbf{D^{b}_{pur}}(R)$, then $X$ admits pure
projective (resp. injective) resolutions.

Now we may define a functor
$$\mathbf{R}\Hom_{R}(-,-):\mathbf{D^{b}_{pur}}(R)^{op}\times
\mathbf{D^{b}_{pur}}(R)\rightarrow
\mathbf{D_{pur}}(\mathbb{Z})$$ using either the pure projective
resolution of the first variable or the pure injective resolution of
the second variable. More precisely, let $P_X$ be a pure projective resolution of $X$ and $I_Y$ a pure injective resolution of $Y$.
Then we have a diagram of pure quasi-isomorphisms
$$\mathbf{R}\Hom_R(X,Y):=\Hom_R(P_X,Y)\rightarrow \Hom_R(P_X,I_Y)\leftarrow \Hom_R(X,I_Y):=\mathbf{R}\Hom_R(X,Y).$$
It follows that $\mathbf{R}\Hom_{R}(-,-)$ is well defined, and we call it the {\it right pure derived functor} of Hom.

Let $P\rightarrow X$ be a pure projective resolution of $X$ and
$Y\rightarrow I$ a pure injective resolution of $Y$. In order to
coincide with the classical ones in [KS] and [S], we put
$$\Pext_{R}^{i}(X,Y):=\H ^{i}\mathbf{R}\Hom_{R}(X,Y)=\H ^{i}\Hom_{R}(P,Y),\ {\rm and}$$
$$\Pext_{R}^{i}(X,Y):=\H ^{i}\mathbf{R}\Hom_{R}(X,Y)=\H ^{i}\Hom_{R}(X,I).$$

Recall that $X\in\mathbf{C}(R)$ is called {\it contractible} if it is isomorphic to the zero object in $\mathbf{K}(R)$,
equivalently, the identical map $\Id_X$ is homotopic to zero. That is to say, $X$ is splitting exact (see [We, Exercise 1.4.3]).

\begin{theorem}
For any $X\in \mathbf{D^{b}_{pur}}(R)$ and $n\in \mathbb{Z}$, the
following statements are equivalent.
\begin{enumerate}
\item[(1)] $\ppd_{R}X\leq n$.

\item[(2)] $\mathbf{inf_{p}}X \geq -n$, and if $f':P'\rightarrow X$ is a pure
projective resolution of $X$, then the R-module $\Coker
d_{P'}^{-n-1}$ is pure projective.

\item[(3)] If $f':P'\rightarrow X$ is a pure projective resolution of $X$,
then $P'=P_{1}\bigoplus P_{2}$, where $P_{1}^{i}=0$ for any $i<-n$ and $P_{2}$ is
contractible.

\item[(4)] $\Pext_{R}^{i}(X,Y)=0$ for any $Y\in \mathbf{D_{pur}}(R)$ and $i>n+
\mathbf{sup_{p}}Y$.

\item[(5)] $\mathbf{inf_{p}}X\geq -n$ and $\Pext_{R}^{n+1}(X,N)=0$ for any $N\in R$-$\Mod$.
\end{enumerate}
\end{theorem}

{\it Proof.} $(1)\Rightarrow (2)$ Let $\ppd_{R}X\leq n$. Then there
exists a pure projective resolution $f:P\rightarrow X$ with
$P^{i}=0$ for any $i<-n$. By Proposition 3.4, we have $\mathbf{inf_{p}}X \geq -n$.
Let $f':P'\rightarrow X$ be another pure projective resolution of
$X$. Then there exists a quasi-isomorphism of complexes
$$\Hom_{R}(P,f'):\Hom _{R}(P,P')\rightarrow\Hom _{R}(P,X).$$ Thus
there exists a cochain map $g: P\to P'$ such that $f'\circ g=f$, and therefore
$$\Hom_R(F,f')\circ \Hom_R(F,g)=\Hom_R(F,f)$$ for any finitely presented $R$-module $F$.
It follows from Remark 2.8(2) that $g$ is a pure quasi-isomorphism. Then $g$ is a homotopy
equivalence by Lemma 2.10. It is easy to check that the exact sequence
$$\CD
\cdots @> >> P'^{-n-1} @>d_{P'}^{-n-1} >> P'^{-n} @>
\widetilde{d}_{P'}^{-n-1}>>\Coker d_{P'}^{-n-1} @> >> 0
\endCD$$ is contractible. So $\Coker d_{P'}^{-n-1}$ is pure projective.

$(2)\Rightarrow (3)$ Let $f':P'\rightarrow X$ be a pure projective resolution of $X$.
Because $\mathbf{inf_{P}}P'=\mathbf{inf_{p}}X \geq -n$, we have that the sequence
$$\CD
\cdots @> >> P'^{-n-1} @> d_{P'}^{-n-1} >> P'^{-n} @> \widetilde{d}_{P'}^{-n-1} >>\Coker d_{P'}^{-n-1} @> >>
0~~~~~~~~~~~~~~(4.1)
\endCD$$ is pure exact. Because $\Coker d_{P'}^{-n-1}$ is pure projective by assumption, $(4.1)$ is
contractible. Now let $P'^{-n}=M\oplus \Coker d_{P'}^{-n-1}$, and put
$$P_1:=~~~~\cdots \to 0\to \Coker d_{P'}^{-n-1}\to P'^{-n+1}\to P'^{-n+2}\to \cdots,\ {\rm and}$$
$$P_2:=~~~~\cdots \to P'^{-n-2} \to P'^{-n-1}\to M\to 0\to \cdots.$$ Then we have $P'=P_{1}\bigoplus
P_{2}$, where $P_{1}^{i}=0$ for any $i<-n$ and $P_{2}$ is contractible.

$(3)\Rightarrow (1)$ By (3), we have that the embedding $P_{1}\hookrightarrow P'$ is
clearly a pure quasi-isomorphism. This implies that $X$ admits a pure
projective resolution $P_{1}\hookrightarrow P'\rightarrow X$ with
$P_{1}^{i}=0$ for any $i<-n$.

$(3)\Rightarrow (4)$ We only need to consider the situation when $\mathbf{sup_{p}}Y=m<\infty$.
Let $P\to X$ be a pure projective resolution of
$X$. Then $P=P_{1}\bigoplus P_{2}$, where $P_{1}^{i}=0$ for any $i<-n$
and $P_{2}$ is contractible. So we have
$$\Pext_{R}^{i}(X,Y)=\H^{i}\Hom_{R}(P,Y)=\H^{i}\Hom_{R}(P_{1},Y).$$
As in Lemma 2.11, let
$Y'$ be the right canonical truncation complex of $Y$ at degree $m$. Then the embedding $Y'\hookrightarrow Y$ is a
pure quasi-isomorphism. So we have
$$\H^{i}\Hom_{R}(P_{1},Y)=\H^{i}\Hom_{R}(P_{1},Y')=0$$ for any $i>n+
m$. Thus $\Pext _{R}^{i}(X,Y)=0$ for any $Y \in
\mathbf{D_{pur}}(R)$ and $i>n+\mathbf{sup_{p}}Y$.

$(4)\Rightarrow (5)$ For any $N\in R$-$\Mod$, we have $\mathbf{sup_{p}}N=0$ and
$\Pext _{R}^{n+1}(X,N)=0$ by (4).

Let $M$ be a right $R$-module. Then $$\H^i((M\otimes_R X)^+)=\H^i(\Hom_R(X,M^+))=\Pext_{R}^{i}(X,M^+)=0$$
for any $i>n$ by the adjoint isomorphism theorem and (4). So $M\otimes_R X$ is exact in degree $<-n$,
and hence $X$ is pure exact in degree $<-n$ by Proposition 3.4. It implies that $\mathbf{inf_{p}}X\geq -n$.

$(5)\Rightarrow (3)$ Let $P'$ be a pure projective resolution of $X$
and $N\in R$-$\Mod$. Then we have
$\mathbf{inf_{p}}P'=\mathbf{inf_{p}}X\geq -n$. So $P'$ is pure exact
in degree $\leq -n-1$, and hence the sequence
$$\CD \cdots @>
>> P'^{-n-2} @> >> P'^{-n-1} @> >> P'^{-n} @> >> \mathrm{Coker}_{P'}^{-n-1} @> >> 0
\endCD$$is pure exact and it is a pure projective resolution of $\mathrm{Coker}_{P'}^{-n-1}$.
We have the following equalities
$$\Pext_{R}^{1}(\mathrm{Coker}_{P'}^{-n-1},N)=\H ^{n+1}\Hom_{R}(P',N)=\Pext^{n+1}_{R}(X,N)=0.$$
It implies that $\mathrm{Coker}_{P'}^{-n-1}$ is pure projective. Thus the above pure
exact complex is contractible, and therefore $P'=P_{1}\bigoplus P_{2}$,
where $P_{1}^{i}=0$ for any $i<-n$. \hfill$\square$

\vspace{0.2cm}

Dually, we have the following

\begin{theorem}
For any $Y\in\mathbf{ D^{b}_{pur}}(R)$ and $n\in \mathbb{Z}$, the
following statements are equivalent.
\begin{enumerate}
\item[(1)] $\pid_{R}Y\leq n$.

\item[(2)] $\mathbf{sup_{p}}Y\leq n$, and if $f':Y\rightarrow I'$ is a pure
injective resolution of $X$, then the R-module $\Ker d_{I'}^{n}$ is
pure injective.

\item[(3)] If $f':Y\rightarrow I'$ is a pure injective resolution of $X$,
then $I'=I_{1}\bigoplus I_{2}$, where $I_{1}^{i}=0$ for any $i>n$ and
$I_{2}$ is contractible.

\item[(4)] $\Pext_{R}^{i}(X,Y)=0$ for any $X\in \mathbf{D_{pur}}(R)$ and $i>n -\mathbf{inf_{p}}X$.

\item[(5)] $\mathbf{sup_{p}}Y\leq n$ and $\Pext_{R}^{n+1}(M,Y)=0$ for any $M\in R$-$\Mod$.
\end{enumerate}
\end{theorem}

By the above two theorems and Proposition 3.4, for any complex $X \in\mathbf{
D^{b}_{pur}}(R)$, we have the following characterizations of $\ppd
_{R}X$ and $\pid _{R}X$ via the pure derived functor $\mathbf{R}\Hom$

$$~~~~~~\ppd _{R}X=\sup\{i\in \mathbb{Z} \mid \Pext _{R}^{i}(X,N)\neq 0
{\rm \ for \ some}\ N\ \textrm{in}\ R{\text-}\Mod\}, \ {\rm and} $$

$$\pid _{R}Y=\sup\{i\in \mathbb{Z} \mid \Pext _{R}^{i}(M,Y)\neq 0
{\rm \ for \ some}\ M\ \textrm{in}\ R{\text-}\Mod\}.$$

\vspace{0.2cm}

Recall that the {\it left pure global dimension} of $R$, written $\Pgldim
R$, is the supremum of the pure projective dimension of all modules
in $R$-$\Mod$. It is also equals to the supremum of the
pure injective dimension of all modules in $R$-$\Mod$. It
is well known that $\Pgldim R\leq n$ if and only if
$\Pext_{R}^{i}(M,N)=0$ for any $M,N\in R$-$\Mod$ and $i>n$, see for example [S, p.95].
We have the cochain complex version of this result.

\begin{theorem} For any $n\in \mathbb{Z}$, the following statements are equivalent.
\begin{enumerate}
\item[(1)] $\Pgldim R\leq n$.

\item[(2)] $\ppd_{R}X\leq n-\mathbf{inf_{p}}X$ for any $X\in
\mathbf{D_{pur}^{b}}(R)$.

\item[(3)] $\pid_{R}Y\leq n+\mathbf{sup_{p}}Y$ for any $Y\in
\mathbf{D_{pur}^{b}}(R)$.

\item[(4)] $\Pext_{R}^{i}(X,Y)=0$ for any $X,Y\in \mathbf{D_{pur}^{b}}(R)$
and $i>n+\mathbf{sup_{p}}Y-\mathbf{inf_{p}}X$.
\end{enumerate}
\end{theorem}

{\it Proof.} The implications $(2)\Rightarrow (4)$ and $(3)\Rightarrow (4)$ follow
from Theorems 4.5 and 4.6, respectively. The implication $(4)\Rightarrow (1)$ is
obvious just letting both $X$ and $Y$ be $R$-modules.
The implication $(1)\Rightarrow (2)$ is the dual of $(1)\Rightarrow (3)$. So it remains to
prove the implication $(1)\Rightarrow (3)$.

Let $\mathbf{sup_{p}}Y=m$ and $Y\rightarrow I$ be a pure injective
resolution of $Y$. Then $\mathbf{sup_{p}}I=m$ and $I$ is pure exact in degree $\geq m+1$.
So $$0 \to \Ker d_{I}^{m}\to I^{m}\to I^{m+1}\to \cdots$$
is a pure injective resolution of $\Ker d_{I}^{m}$. By (1), we have $\pid_{R}\Ker d_{I}^{m}\leq n$. Let
$$0\rightarrow \Ker d_{I}^{m}\to K^{0}\to K^{1} \to \cdots \to K^{n}\to 0$$
be a pure injective resolution of $\Ker d_{I}^{m}$. Then it is easy to check that
$$\cdots \to I^{m-2} \to I^{m-1} \to K^{0} \to \cdots \to K^{n-1} \to
 K^{n} \to 0 \to\cdots$$ is a pure injective resolution of $Y$ and $\pid_{R}Y\leq n+m$. \hfill$\square$

\bigskip
\section{ The case of unbounded complexes}
\setcounter{equation}{0}

In this section, we study the existence of pure projective resolutions of unbounded
complexes. We need the tools of homotopy colimits and limits ([N2]).

Let
$$\CD X_{0} @>j_{1}>> X_{1} @>j_{2}>> X_{2} @>j_{3}>> \cdots ~~~~~~~~~~~~~~~~~~~~~~~~~~~~~~~~~~~~~~~~(5.1)\endCD$$
be a sequence in $\mathbf{K}^b(R)$, where $j_{i}$
is a morphism of complexes for any $i>0$. Then we can form the {\it
homotopy colimit} of this sequence, written
$\Ho\underrightarrow{\colim}X_{i}$, by the triangle
$$\CD
\bigoplus\limits_{i=0}^{\infty} X_{i} @> \textrm{1-shift} >> \bigoplus\limits_{i=0}^{\infty} X_{i} @> >>
\Ho \underrightarrow{\colim}X_{i} @> >> (\bigoplus\limits_{i=0}^{\infty}
X_{i})[1]
\endCD$$ in $\mathbf{K}(R)$. The notion of {\it homotopy limits} is defined dually, and denoted by $\Ho\underleftarrow{\mathrm{lim}}$.

For the sequence (5.1), we can also form the {\it direct limit}, written $\underrightarrow{\colim} X_{i}$, in
$\mathbf{C}(R)$. We have the following exact sequence of complexes
(note: the morphism $\textrm{1-shift}$ is monic)
$$\CD 0 @> >> \bigoplus\limits_{i=0}^{\infty} X_{i} @>\textrm{1-shift}>>
\bigoplus\limits_{i=0}^{\infty} X_{i} @> \iota >> \underrightarrow{\colim} X_{i} @>
>>0.~~~~~~~~~~~~~~(5.2) \endCD$$
Then it is easy to check that there exists a morphism $\alpha:
\Ho\underrightarrow{\colim}X_{i}\to\underrightarrow{\colim} X_{i}$, since $\iota
\circ(\textrm{1-shift})=0$. Passing to homology we conclude that $\alpha$
is a quasi-isomorphism, since H$^0(-)$ is cohomological by [GM, Chapter IV.1.6].

Because $\Hom_R(F,-)$ commutes with direct limits in $\mathbf{C}(R)$ for any finitely presented $R$-module $F$ by
[AR, Corollary 1.54] (see also [St, Remark 4.13] or [CH, Corollary 4.6]), we have that (5.2) is pure exact in $\mathbf{C}(R)$,
that is, it is pure exact in each degree. So after
applying the functor $\Hom_{R}(P,-)$ for any $P\in \mathcal
{P}\mathcal {P}$, we get the following exact sequence
$$\CD 0 @> >> \Hom_{R}(P,\bigoplus\limits_{i=0}^{\infty} X_{i}) @>\Hom_{R}(P,\textrm{1-shift})>>
\Hom_{R}(P,\bigoplus\limits_{i=0}^{\infty} X_{i} )@> >> \Hom_{R}(P,\underrightarrow{\colim} X_{i}) @>
  >>0
\endCD,$$
and the following exact triangle
{\footnotesize $$\CD
\Hom_{R}(P, \bigoplus\limits_{i=0}^{\infty} X_{i}) @> \Hom_{R}(P,\textrm{1-shift}) >>
\Hom_{R}(P,\bigoplus\limits_{i=0}^{\infty} X_{i}) @> >>\Hom_{R}(P, \Ho
\underrightarrow{\colim}X_{i})  @> >>
\Hom_{R}(P,\bigoplus\limits_{i=0}^{\infty} X_{i})[1]
\endCD$$} in $\mathbf{K}(R)$. When passing to homology we see that $\alpha$
is a pure quasi-isomorphism.

\begin{theorem}
Let $X\in\mathbf{C}(R)$ be a bounded below complex. Then there exists a complex $P$
consisting of pure projective $R$-modules satisfying the following properties.
\begin{enumerate}
\item[(1)] There exists a pure quasi-isomorphism $f:P\rightarrow X$.

\item[(2)] $\Hom_R(P,-)$ preserves pure exact complexes.
\end{enumerate}
That is, $f:P\to X$ is a pure projective resolution of $X$.
\end{theorem}

{\it Proof.} (1) Write $X:=\underrightarrow{\colim} X_{i}$ with the structure map $j_{i+1}:X_{i}\to X_{i+1}$, where $X_{i}$ is a bounded
complex for any $i\geq 0$. By Proposition 3.5, for any $X_{i}$ there
exists a pure quasi-isomorphism $f_{i}:P_{i}\longrightarrow X_{i}$
with $P_{i}\in \mathbf{K}^{-}(\mathcal{PP})$ for any $i\geq 0$. Then we obtain the following commutative diagram
$$\CD
P_{i} @>\overline{j_{i+1}}>> P_{i+1} \\
@V f_{i} VV @V f_{i+1} VV  \\
X_{i} @>j_{i+1}>> X_{i+1}
\endCD$$
in $\mathbf{K}(R)$, where $\overline{j_{i+1}}$ is induced by
$j_{i+1}$. So there exists a morphism of exact triangles
$$\CD
\bigoplus\limits_{i=0}^{\infty} P_{i} @>\textrm{1-shift}>> \bigoplus\limits_{i=0}^{\infty} P_{i} @> >>
\Ho\underrightarrow{\colim} P_{i} @> >> (\bigoplus\limits_{i=0}^{\infty} P_{i})[1] \\
@V    VV @V    VV @V f  VV @V   VV  \\
\bigoplus\limits_{i=0}^{\infty} X_{i} @>\textrm{1-shift}>> \bigoplus\limits_{i=0}^{\infty} X_{i} @> >>
\Ho\underrightarrow{\colim} X_{i} @> >> (\bigoplus\limits_{i=0}^{\infty} X_{i})[1]
\endCD$$ in
$\mathbf{K}(R)$. After applying the localization functor, it is a morphism of exact
triangles in $\mathbf{D_{pur}}(R)$. Since pure
quasi-isomorphisms are closed under coproducts by Remark 2.8(2),
we have that the first two vertical maps in the above diagram are pure quasi-isomorphisms.
So $f$ and
$$\alpha\circ f:P=\Ho\underrightarrow{\colim} P_{i}\rightarrow
\Ho\underrightarrow{\colim} X_{i}\rightarrow \underrightarrow{\colim} X_{i}$$ 
are also pure quasi-isomorphisms. By the
construction, we have that $\Ho\underrightarrow{\colim} P_{i}$ is the mapping
cone of some cochain map between complexes consisting of pure
projective $R$-modules. Thus $\Ho\underrightarrow{\colim} P_{i}$ is also
a complex consisting of pure projective $R$-modules.

(2) We will prove that $\Hom_{\mathbb{Z}}(F,\Hom_{R}(P,X))$ is exact
for any pure exact complex $X$ of $R$-modules and any finitely
presented $\mathbb{Z}$-module $F$. Consider the following commutative diagram
$$\CD
\Hom_{R}(P,X) @> >> \Hom_{R}(\bigoplus\limits_{i=0}^{\infty} P_{i},X) @> >>
\Hom_{R}(\bigoplus\limits_{i=0}^{\infty} P_{i},X) @> >> (\Hom_{R}(P,X))[1] \\
@V =  VV @V  \cong VV @V \cong VV @V   VV  \\
\Hom_{R}( P,X) @> >> \prod\limits_{i=0}^{\infty} \Hom_{R}( P_{i},X) @> >> \prod\limits_{i=0}^{\infty} \Hom_{R}(
P_{i},X) @> >> (\Hom_{R}(
P,X))[1] \endCD$$ in $\mathbf{K}(\mathbb{Z})$, where both rows are
exact triangles. We have the following isomorphisms
$$\Hom_{\mathbb{Z}}(F,\Hom_{R}(\bigoplus\limits_{i=0}^{\infty}
P_{i},X))\cong \prod \limits_{i=0}^{\infty}
\Hom_{\mathbb{Z}}(F,\Hom_{R}(P_{i},X)).$$ Note that
the latter one is exact by Remark 4.2. Because $\Hom_{\mathbb{Z}}(F,-)$ is a triangulated functor,
the assertion follows standardly.
\hfill$\square$

\begin{theorem}
Let $X\in\mathbf{C}(R)$ be a bounded above complex. Then there exists a complex $I$
consisting of pure injective $R$-modules satisfying the following properties.
\begin{enumerate}
\item[(1)] There exists a pure quasi-isomorphism $f:X\rightarrow I$.

\item[(2)] $\Hom_R(-,I)$ preserves pure exact complexes.
\end{enumerate}
That is, $f:X\to I$ is a pure injective resolution of $X$.
\end{theorem}

{\it Proof.} Write $X:=\underleftarrow{\mathrm{lim}} X_{i}$ with $X_{i}$ a bounded complex for any $i\leq 0$.
Then by [IK, Lemma 2.6], we have $X\cong \Ho\underleftarrow{\mathrm{lim}} X_{i}$
in $\mathbf{K}(R)$. Note that pure quasi-isomorphisms are closed under products by Remark 2.8(2).
Now by using an argument similar to that in the proof of Theorem 5.1, we get the assertion. \hfill$\square$

\begin{remark}{\rm
One can find that the derived functor
$$\mathbf{R}\Hom_{R}(-,-):\mathbf{D^{b}_{pur}}(R)^{op}\times
\mathbf{D^{b}_{pur}}(R)\rightarrow
\mathbf{D_{pur}}(\mathbb{Z})$$
may be extended to
$$\mathbf{R}\Hom_{R}(-,-):\mathbf{D^{+}_{pur}}(R)^{op}\times
\mathbf{D^{-}_{pur}}(R)\rightarrow
\mathbf{D_{pur}}(\mathbb{Z}).$$ The
corresponding characterizations of dimensions in Section 4 also hold in this situation.}
\end{remark}

\vspace{0.5cm}
\noindent $\mathbf{Acknowledgements}$

This research was partially supported by NSFC (Grant Nos. 11171142 and 11571164)
and a Project Funded by the
Priority Academic Program Development of Jiangsu Higher Education Institutions.
The authors thank the referee for the helpful suggestions.

\end{document}